\documentclass[11pt]{article}
\usepackage{amsmath, epsfig, cite, amssymb, amsfonts,latexsym}
\newcommand{\qed}{{\hfill\rule{4pt}{7pt}}}

\newtheorem{thm}{Theorem}[section]

\newtheorem{pro}[thm]{Proposition}

\def\pf{\noindent {\it Proof.} }

\makeatletter \@addtoreset{equation}{section} \makeatother

\begin{document}
\begin{center}
{\bf \Large Semi-Finite Forms of

Bilateral Basic Hypergeometric Series}
\end{center}
\begin{center}
{ William Y. C. Chen}$^{1}$ \quad  and  \quad { Amy M. Fu}$^{2}$

   Center for Combinatorics, LPMC\\
   Nankai University, Tianjin 300071, P.R. China

   \vskip 1mm

   Email: $^1$chen@nankai.edu.cn, $^2$fu@nankai.edu.cn
\end{center}

\noindent{\bf Abstract.} We show that several classical bilateral
summation and transformation formulas have semi-finite forms. We
obtain these semi-finite forms from unilateral summation and
transformation formulas. Our method can be applied to derive
Ramanujan's $_1\psi_1$ summation, Bailey's $_2\psi_2$
transformations, and Bailey's $_6\psi_6$ summation.

\noindent {\bf Corresponding Author:} William Y. C. Chen, Email:
chen@nankai.edu.cn



 \vskip 3mm

\section{Introduction}

We follow the terminology for basic hypergeometric series in
\cite{GR}. Assuming $|q|<1$, let
\[ (a;q)_\infty = (1-a) (1-aq) (1-aq^2) \cdots .\]
For any integer $n$, the $q$-shifted factorial $(a;q)_n$ is given
by
\[ (a;q)_n = { (a;q)_\infty \over (aq^n;q)_\infty}.\]
For $n\geq 0$, we have the following relation which is crucial for
this paper:
\begin{equation}\label{Defi}
(a;q)_{-n}= \frac{1}{(aq^{-n};q)_n}={ (-q/a)^{n} q^{\binom{n}{2}
}\over (q/a;q)_{n}} .
\end{equation}
For convenience, we employ the following usual notation:
$$
(a_1, a_2, \ldots, a_m;q)_n=(a_1;q)_n(a_2;q)_n\ldots(a_m;q)_n.
$$
 The (unilateral) basic
hypergeometric series $_{r+1}\phi_r$ is defined by
\begin{eqnarray}\label{Hype}
_{r+1}\phi_r\left[
\begin{array}{c}
a_1, a_2, \cdots, a_{r+1}\\
b_1, b_2, \cdots, b_{r}
\end{array};q, z
\right]=\sum_{k=0}^{\infty}A(k),
\end{eqnarray}
where
$$
A(k)=\frac{(a_1, a_2, \cdots, a_{r+1};q)_k}{(b_1, b_2, \cdots,
b_r,q;q)_k}z^k.
$$
 The bilateral basic hypergeometric series $_s\psi_s$ is
defined as follows,
\begin{eqnarray}\label{Bila}
_s\psi_s\left[
\begin{array}{c}
a_1, a_2, \cdots, a_{s}\\
b_1, b_2, \cdots, b_{s}
\end{array};q, z
\right]=\sum_{k=-\infty}^{\infty} B(k),
\end{eqnarray}
where
$$
B(k)=\frac{(a_1, a_2, \cdots, a_{s};q)_k}{(b_1, b_2, \cdots,
b_s;q)_k}z^k.
$$

In this paper, we propose the following method of deriving
bilateral summation and transformation formulas using {\em
semi-finite forms}. For a  bilateral series $_s\psi_s$ as given in
(\ref{Bila}), we construct a summand $G(k,m)$ which implies a
unilateral series $_{r+s+1}\phi_{r+s}$, where $r$ is a nonnegative
integer, such that
\[ \lim _{m \rightarrow \infty }G(k,m)=B(k)\]
for all $k$, and the summation
\begin{equation}\label{gn0}
 \sum_{k=-m}^\infty G(k,m)\end{equation}
can be easily accomplished  as a Laurent extension of the
summation
\begin{equation}\label{laurants}
\sum_{k=0}^\infty  G(k-m, m)= G(-m,m)\sum_{k=0}^\infty A(k),
\end{equation} where $G(k,m)$ can be written as
 \[ G(k-m, m) = G(-m,m) A(k)\]
 for some $A(k)$. The bilateral series (\ref{Bila}) is then obtained
    from (\ref{gn0}) as $m\to\infty$, subject to suitable convergence
    conditions. We apply this procedure  to derive
    bilateral series identities from suitable unilateral ones.
    The above summation (\ref{gn0}) is called  the {\em
    semi-finite form} of the bilateral summation (\ref{Bila}).
  A method similar to ours was recently used by Schlosser \cite{SCHL},
    and Jouhet and Schlosser \cite{SCHL04}, who derived summations for bilateral
    series from {\em finite forms}.
    We also note that another method, which uses a similar
    factorization as above, for deriving bilateral series identities
     from unilateral ones was used by Ismail \cite{Ismail}, and Askey and
 Ismail \cite{AsIs}. Rather than taking limits, they
 apply analytic continuation as the main ingredient.

 In this paper, we present semi-finite forms of several classical bilateral summation
and transformation formulas such as Ramanujan's $_{1}\psi_1$
formula, Bailey's $_2\psi_2$ transformations,  and Bailey's
$_6\psi_6$ summation.

\section{From $_2\phi_1$ to $_1\psi_1$}
Using the well known Gauss summation formula
\begin{equation}\label{Gauss}
_2\phi_1\left[
\begin{array}{c}
a,b\\
c
\end{array};q, c/ab
\right]=\frac{(c/a,c/b;q)_{\infty}}{(c,c/ab;q)_{\infty}},
\end{equation}
where $|c/ab|<1$, we get a semi-finite form of Ramanujan's
summation of the general $_1\psi_1$,
\begin{equation}\label{Ran}
_{1}\psi_1 \left[
\begin{array}{l}
a\\
b
\end{array};q,z\right]=\sum_{k=-\infty}^{\infty}\frac{(a;q)_k}{(b;q)_k}z^k=
\frac{(q;q)_{\infty}(b/a;q)_{\infty}(az;q)_{\infty}(q/az;q)_{\infty}}
{(b;q)_{\infty}(q/a;q)_{\infty}(z;q)_{\infty}(b/az;q)_{\infty}},
\end{equation}
where $|b/a|<|z|<1$.

\begin{pro} For $|z|<1$, the following identity holds:
 \label{theo}
\begin{equation} \label{r-f}
\sum_{k=-m}^{\infty}\frac{(a;q)_k(bq^{m}/az;q)_k}{(q^{1+m};q)_k(b;q)_k}z^k
=\frac{(q;q)_m(q/az;q)_m}{(q/a;q)_m(b/az;q)_m}
\frac{(b/a;q)_{\infty}(az;q)_{\infty}}{(b;q)_{\infty}(z;q)_{\infty}}.
\end{equation}
\end{pro}
\noindent{\it Proof.} The left hand side of (\ref{r-f}) can be
rewritten as
 \begin{eqnarray*}
\lefteqn{\sum_{k=0}^{\infty}\frac{(a;q)_{k-m}(bq^{m}/az;q)_{k-m}}
{(q^{1+m};q)_{k-m}(b;q)_{k-m}}z^{k-m}}\\[6pt]
&=&z^{-m}\frac{(a;q)_{-m}(bq^{m}/az;q)_{-m}}{(q^{1+m};q)_{-m}(b;q)_{-m}}
\sum_{k=0}^{\infty}\frac{(aq^{-m};q)_k(b/az;q)_k}{(q;q)_k(bq^{-m};q)_k}z^k\\[6pt]
&\overset{(\ref{Gauss})}{=}&z^{-m}\frac{(a;q)_{-m}(bq^{m}/az;q)_{-m}}{(q^{1+m};q)_{-m}(b;q)_{-m}}
\frac{(b/a;q)_{\infty}(azq^{-m};q)_{\infty}}{(bq^{-m};q)_{\infty}(z;q)_{\infty}}
 \\[6pt]
&\overset{(\ref{Defi})}{=}&z^{-m}\frac{(q;q)_{m}(azq^{-m};q)_m}{(aq^{-m};q)_m
(b/az;q)_m}\frac{(az;q)_{\infty}(b/a;q)_{\infty}}
{(b;q)_{\infty}(z;q)_{\infty}},
\end{eqnarray*}
which equals the right hand side of (\ref{r-f}). \qed

Taking the limit $m  \rightarrow \infty$ in Proposition \ref{theo}
while assuming $|b/az|<1$, we immediately obtain (\ref{Ran}).

We remark that our method is different from the method of M.
Jackson's elementary proof of (\ref{Ran}) (see the exposition of
Schlosser \cite{SCHL}) in the sense that Jackson's proof does not
give a semi-finite form although the Gauss summation is also the
basic ingredient.  We should also note that a finite form of
Ramanujan's $_1\psi_1$ summation has been given by Schlosser
\cite{Schl03b} using the terminating $q$-Pfaff-Saalsch\"utz
summation.

\section{From $_3\phi_2$ to $_2\psi_2$}

In this section, we  use two $_3\phi_2$ summation and
transformation formulas to give  the semi-finite forms of
$_2\psi_2$ formulas due to Bailey. We begin with the following
$_2\psi_2$ transformation formula \cite[Ex. 5.20(i)]{GR} valid for
$|z|, |cd/abz|,|d/a|,|c/b|<1$:
\begin{eqnarray}\label{22}
_{2}\psi_{2}\left[
\begin{array}{l}
a,b\\
c,d
\end{array};q,z
\right]=\frac{(az,d/a,c/b,dq/abz;q)_{\infty}}{(z,d,q/b,cd/abz;q)_{\infty}}
\, {_{2}}\psi_{2}\left[
\begin{array}{l}
a,abz/d\\
az,c
\end{array};q,\frac{d}{a}
\right].
\end{eqnarray}

Using a $q$-analogue of the Kummer-Thomae-Whipple formula
\cite[Eq. (3.2.7)]{GR}:
\begin{eqnarray}\label{Kummer}
_{3}\phi_2\left[
\begin{array}{l}
a, b, c\\
d,e
\end{array};q, \frac{de}{abc}
 \right]=\frac{(e/a,de/bc;q)_{\infty}}{(e,de/abc;q)_{\infty}}{_{3}}\phi_2\left[
\begin{array}{l}
a, d/b, d/c\\
d,de/bc
\end{array};q, \frac{e}{a}
 \right],
\end{eqnarray}
where $|de/abc|<1$ and $|e/a|<1$, we get a semi-finite form of
(\ref{22}).

\begin{pro}\label{Thm1} For $|z|<1$ and $|d/a|<1$, we have
 \begin{eqnarray}
\sum_{k=-m}^{\infty}\frac{(a,b;q)_k(cdq^m/abz;q)_k}{(c,d;q)_k(q^{1+m};q)_k}z^k&=&
\frac{(az,d/a;q)_{\infty}}{(z,d;q)_{\infty}}\frac{(c/b,dq/abz;q)_{m}}
{(q/b,cd/abz;q)_{m}}  \nonumber  \\\label{xy}
 &&
    \;\; \; \cdot \sum_{k=-m}^{\infty}\frac{(a,cq^m/b,abz/d;q)_{k}}
    {(c,q^{1+m},az;q)_{k}}(d/a)^{k}.
\end{eqnarray}
\end{pro}

\pf  The left hand side of (\ref{xy}) equals
\begin{eqnarray*}
\lefteqn{z^{-m}\frac{(a,b,cdq^m/abz;q)_{-m}}{(c,d,q^{1+m};q)_{-m}}\sum_{k=0}^{\infty}\frac{(aq^{-m},bq^{-m},cd/abz;q)_k}
{(cq^{-m},dq^{-m},q;q)_k}z^k}
\\
&\overset{(\ref{Kummer})}{=}&z^{-m}\frac{(a,b,cdq^m/abz;q)_{-m}}{(c,d,q^{1+m};q)_{-m}}\frac{(d/a,azq^{-m};q)_{\infty}}{(dq^{-m},z;q)_{\infty}}\\
&& \;\;\; \cdot \sum_{k=0}^{\infty}\frac{(aq^{-m},c/b,
abzq^{-m}/d;q)_k}{(q,cq^{-m},azq^{-m};q)_k}\left(\frac{d}{a}\right)^k\\
&\overset{(\ref{Defi})}{=}&\frac{(d/a,az;q)_{\infty}}{(d,z;q)_{\infty}}\frac{(c/b,abzq^{-m}/d;q)_m}
{(bq^{-m}, cd/abz;q)_m}\left(\frac{d}{az}\right)^m\\
& & \;\;\; \cdot
\sum_{k=0}^{\infty}\frac{(a,cq^m/b,abz/d;q)_{k-m}}
{(c,q^{1+m},az;q)_{k-m}}(d/a)^{k-m},
\end{eqnarray*}
which can be rewritten in the form of the right hand side of
(\ref{xy}).  \qed

The next $_2\psi_2$ transformation formula we consider is the
following \cite[Ex. 5.20(ii)]{GR}:
\begin{eqnarray} \label{b2}
_{2}\psi_{2}\left[
\begin{array}{c}
a,b\\
c,d
\end{array};q,z
\right]=\frac{(az,bz,cq/abz,dq/abz;q)_{\infty}}{(q/a,q/b,c,d;q)_{\infty}}
\,  {_{2}}\psi_{2}\left[
\begin{array}{l}
abz/c,abz/d\\
az,bz
\end{array};q,\frac{cd}{abz}
\right] .
\end{eqnarray}

Using a summation of Hall \cite[Eq. (3.2.10)]{GR}:
\begin{eqnarray}\label{Hall}
_{3}\phi_2\left[
\begin{array}{l}
a, b, c\\
d,e
\end{array};q, \frac{de}{abc}
 \right]=\frac{(b,de/ab,de/bc;q)_{\infty}}{(d,e,de/abc;q)_{\infty}}{_{3}}\phi_2\left[
\begin{array}{l}
d/b,e/b, de/abc\\
de/ab,de/bc
\end{array};q, b
 \right],
\end{eqnarray}
where $|de/abc|<1$ and $|b|<1$, we obtain the following
semi-finite form of (\ref{b2}).
\begin{pro} For $|z|<1$ and $|cd/abz|<1$, we have
\begin{eqnarray*}
\sum_{k=-m}^{\infty}\frac{(a,b;q)_k(cdq^m/abz;q)_k}{(c,d;q)_k(q^{1+m};q)_k}z^k&=&
\frac{(az,bz,cd/abz;q)_{\infty}}{(c,d,z;q)_{\infty}}\frac{(cq/abz,dq/abz,z;q)_m}{(q/a,q/b,cd/abz;q)_m}\\
&&  \;\;\; \cdot
\sum_{k=-m}^{\infty}\frac{(abz/c,abz/d,zq^m;q)_{k}}{(az,bz,q^{1+m};q)_{k}}(cd/abz)^{k}.
\end{eqnarray*}
\end{pro}

\section{From nonterminating $_8\phi_7$ to $_6\psi_6$}

In this section, we give a semi-finite form of Bailey's $_6\psi_6$
summation formula by using Bailey's $3$-term transformation
formula for a nonterminating very-well-poised $_8\phi_7$ series
\cite[Eq. (2.11.1)]{GR}:
\begin{eqnarray}\label{Watson2}
\lefteqn{_{8}\phi_7\left[
\begin{array}{c}
a, qa^{1\over 2},-qa^{1\over 2},b,c, d, e, f\\
a^{1\over 2},-a^{1\over 2},aq/b,aq/c,aq/d,aq/e,aq/f
\end{array};q, \frac{a^2q^2}{bcdef}
 \right] }\nonumber\\
&& \,
=\frac{(aq,aq/de,aq/df,aq/ef,eq/c,fq/c,b/a,bef/a;q)_{\infty}}
{(aq/d,aq/e,aq/f,aq/def,q/c,efq/c,be/a,bf/a;q)_{\infty}}\nonumber\\
&& \qquad \cdot_{8}\phi_7\left[
\begin{array}{c}
ef/c, q(ef/c)^{1\over 2},-q(ef/c)^{1\over 2},aq/bc,aq/cd, ef/a, e, f\\
(ef/c)^{1\over 2},-(ef/c)^{1\over 2}, bef/a,def/a,aq/c,fq/c,eq/c
\end{array};q, \frac{bd}{a}
 \right] \nonumber\\
 &&\qquad +\frac{b}{a}\frac{(aq,bq/a,bq/c,bq/d,bq/e,bq/f,d,e,f;q)_{\infty}}
 {(aq/b,aq/c,aq/d,aq/e,aq/f,bd/a,be/a,bf/a,def/a;q)_{\infty}}\nonumber\\
&&\qquad \cdot \frac{(aq/bc,bdef/a^2,a^2q/bdef;q)_{\infty}}{(aq/def,q/c,b^2q/a;q)_{\infty}}\nonumber\\
&& \qquad \cdot _{8}\phi_7\left[
\begin{array}{c}
b^2/a, qba^{-{1\over 2}},-qba^{-{1\over 2}},b,bc/a, bd/a, be/a, bf/a\\
ba^{-{1\over 2}},-ba^{-{1\over 2}},bq/a,bq/c,bq/d,bq/e,bq/f
\end{array};q, \frac{a^2q^2}{bcdef}
 \right],
\end{eqnarray}
where $|bd/a|<1$ and $|a^2q^2/bcdef|<1$.

\begin{pro} When $|bd/a|<1$ and $|a^2q^2/bcdef|<1$, we have
\begin{eqnarray}\label{66}
\lefteqn{\sum_{k=-m}^{\infty}\frac{(q^{m-k+1}/a,fq^m;q)_{k}}{(q^{m-k}f/a,q^{1+m};q)_k}
\frac{(qa^{1\over 2},-qa^{1\over 2},b,c,d,e;q)_k}{(a^{1\over
2},-a^{1\over 2},aq/b,aq/c,aq/d,aq/e;q)_k}
\left(\frac{qa^2}{bcde}\right)^k}\nonumber\\
&=&\frac{1-efq^{2m}/c}{1-efq^m/c}\frac{(q/a,df/a,ef/a,aq/bc,aq/cd,efq^m/a;q)_m}
{(f/a,q/b,q/c,q/d,def/a,fq^{1+m}/c;q)_m}\nonumber\\
&&
\times\frac{(aq,aq/de,aq/df,aq/ef,eq^{1+m}/c,fq^{1+m}/c,b/a,befq^m/a;q)_{\infty}}
{(aq/d,aq/e,aq/f,aq/def,q^{1+m}/c,efq^{1+m}/c,be/a,bfq^m/a;q)_{\infty}}\nonumber\\
&& \times\sum_{k=-m}^{\infty}\frac{(efq^m/c,q^{1+m}(ef/c)^{1\over
2},-q^{1+m}(ef/c)^{1\over 2},aq^{1+m}/bc;q)_k}
{(q^{1+m},q^{m}(ef/c)^{1\over 2},-q^{m}(ef/c)^{1\over
2},befq^m/a;q)_k}\nonumber\\
&& \qquad
\cdot\frac{(aq^{1+m}/cd,efq^{2m}/a,e,fq^m;q)_k}{(defq^m/a,aq/c,fq^{1+2m}/c,eq^{1+m}/c;q)_k}
\left(\frac{bd}{a}\right)^k\nonumber\\
&&+\frac{b}{a}\frac{1-b^2q^{2m}/a}{1-b^2q^m/a}\left(\frac{a^2q}{bcde}\right)^m
\frac{(q/a,bc/a;q)_{m}}{(f/a;q)_m}\frac{(aq,bq^{1+2m}/a,bq^{1+m}/c;q)_\infty}
{(aq/b,aq/c,aq/d;q)_\infty}\nonumber\\
&& \times
\frac{(bq^{1+m}/d,bq^{1+m}/e,bq/f,d,e,fq^m,aq/bc,bdef/a^2,a^2q/bdef;q)_{\infty}}
{(aq/e,aq/f,bdq^m/a,beq^m/a,bfq^{2m}/a,def/a,aq/def,q/c,b^2q^{1+m}/a;q)_{\infty}}\nonumber\\
&& \times \sum_{k=-m}^{\infty}\frac{(b^2q^m/a,q^{1+m}ba^{-{1\over
2}},-q^{1+m}ba^{-{1\over 2}};q)_k}
{(q^{1+m},q^{m}ba^{-{1 \over 2}},-q^mba^{-{1 \over 2}};q)_k}\nonumber\\
&&\qquad
\cdot\frac{(b,bcq^m/a,bdq^m/a,beq^m/a,bfq^{2m}/a;q)_k}{(bq^{1+2m}/a,bq^{1+m}/c,bq^{1+m}/d,bq^{1+m}/e,bq/f;q)_k}
\left(\frac{a^2q^2}{bcdef}\right)^k.
\end{eqnarray}

\end{pro}

\pf   \, The left hand side of (\ref{66}) equals
\begin{eqnarray*}
\lefteqn{\sum_{k=-m}^{\infty}\frac{(aq^{-m},fq^m;q)_{k}}{(aq^{1-m}/f,q^{1+m};q)_k}
\frac{(qa^{1\over 2},-qa^{1\over 2},b,c,d,e;q)_k}{(a^{1\over
2},-a^{1\over 2},aq/b,aq/c,aq/d,aq/e;q)_k}
\left(\frac{q^2a^2}{bcdef}\right)^k} \\
&\overset{(\ref{laurants})}{=}& \frac{(aq^{-m},fq^m,qa^{1\over
2},-qa^{1\over 2},b,c,d,e;q)_{-m}}{(aq^{1-m}/f,q^{1+m},a^{1\over
2},-a^{1\over 2},aq/b,aq/c,aq/d,aq/e;q)_{-m}}
\left(\frac{q^2a^2}{bcdef}\right)^{-m}\\
&&\times \sum_{k=0}^{\infty}\frac{(aq^{-2m}, q^{1-m}a^{1\over
2},-q^{1-m}a^{1\over 2};q)_{k}}{(q,q^{-m}a^{1\over
2},-q^{-m}a^{1\over
2};q)_k}\\
&&  \qquad \, \cdot
\frac{(bq^{-m},cq^{-m},dq^{-m},eq^{-m},f;q)_k}{(aq^{1-m}/b,aq^{1-m}/c,aq^{1-m}/d,aq^{1-m}/e,aq^{1-2m}/f;q)_k}
\left(\frac{q^2a^2}{bcdef}\right)^k\\
&\overset{(\ref{Watson2})}{=}&\frac{(aq^{-m},fq^m,qa^{1\over
2},-qa^{1\over 2},b,c,d,e;q)_{-m}}{(aq^{1-m}/f,q^{1+m},a^{1\over
2},-a^{1\over 2},aq/b,aq/c,aq/d,aq/e;q)_{-m}}
\left(\frac{q^2a^2}{bcdef}\right)^{-m}\\
&& \times
\frac{(aq^{1-2m},aq/de,aq^{1-m}/df,aq^{1-m}/ef,eq/c,fq^{1+m}/c,bq^m/a,bef/a;q)_{\infty}}
{(aq^{1-m}/d,aq^{1-m}/e,aq^{1-2m}/f,aq/def,q^{1+m}/c,efq/c,be/a,bfq^m/a;q)_{\infty}}\\
&& \times  _{8}\phi_7\left[
\begin{array}{c}
ef/c, q(ef/c)^{1\over 2},-q(ef/c)^{1\over 2},aq/bc,aq/cd, efq^m/a, eq^{-m}, f\\
(ef/c)^{1\over 2},-(ef/c)^{1\over 2},
bef/a,def/a,aq^{1-m}/c,fq^{1+m}/c,eq/c
\end{array};q, \frac{bd}{a}
 \right]\\
 && +\frac{(aq^{-m},fq^m,qa^{1\over 2},-qa^{1\over
2},b,c,d,e;q)_{-m}}{(aq^{1-m}/f,q^{1+m},a^{1\over 2},-a^{1\over
2},aq/b,aq/c,aq/d,aq/e;q)_{-m}}
\left(\frac{q^2a^2}{bcdef}\right)^{-m}\\
&& \times
\frac{bq^m}{a}\frac{(aq^{1-2m},bq^{1+m}/a,bq/c,bq/d,bq/e,bq^{1-m}/f,dq^{-m};q)_{\infty}}
{(aq^{1-m}/b,aq^{1-m}/c,aq^{1-m}/d,aq^{1-m}/e,aq^{1-2m}/f,bd/a,be/a;q)_{\infty}}\\
&& \times
\frac{(eq^{-m},f,aq/bc,bdefq^m/a^2,a^2q^{1-m}/bdef;q)_{\infty}}
{(bfq^m/a,def/a,aq/def,q^{1+m}/c,b^2q/a;q)_\infty}\\
&& \times _{8}\phi_7\left[
\begin{array}{c}
b^2/a, qba^{- {1\over 2}},-qba^{- {1\over 2}},bq^{-m},bc/a,bd/a,be/a, bfq^m/a\\
ba^{-{ 1\over 2}},-ba^{-{1\over 2}},
bq^{1+m}/a,bq/c,bq/d,bq/e,bq^{1-m}/f
\end{array};q, \frac{a^2q^2}{bcdef}
 \right],
\end{eqnarray*}
which equals to the right hand side of (\ref{66}).
 \qed

The above proposition  can be viewed as a semi-finite form of
 Bailey's $_6\psi_6$ summation formula.
 By taking $f=b$ and $m  \rightarrow \infty$ in
(\ref{66}) while assuming $|a^2q/bcde|< 1$, we get

\begin{eqnarray*}
\lefteqn{_{6}\psi_6\left[
\begin{array}{c}
 qa^{1\over 2},-qa^{1\over 2},b,c, d, e\\
a^{1\over 2},-a^{1\over 2},aq/b,aq/c,aq/d,aq/d,aq/e
\end{array};q, \frac{a^2q}{bcde}
 \right] }\\
&=&\frac{(q/a,bd/a,aq/bc,aq/cd,aq, aq/de, aq/bd,
aq/be;q)_{\infty}}{(q/b,q/c,q/d,aq/b,
aq/d,aq/e,aq/bde,bde/a;q)_{\infty}}\\
&& \times
\sum_{k=-\infty}^{\infty}\frac{(e;q)_k}{(aq/c;q)_k}(bd/a)^k\\
&\overset{(\ref{Ran})}{=}&\frac{(q/a,bd/a,aq/bc,aq/cd,aq, aq/de,
aq/bd, aq/be;q)_{\infty}}{(q/b,q/c,q/d,aq/b,
aq/d,aq/e,aq/bde,bde/a;q)_{\infty}}\\
&& \times \frac{(q,aq/ce,bde/a,aq/bde;q)_{\infty}}
{(aq/c,q/e,bd/a,a^2q/bcde;q)_{\infty}}\\
 &=&\frac{(aq,aq/bc,aq/bd,aq/be,aq/ce,aq/cd,aq/de,q,q/a;q)_{\infty}}
 {(aq/b,aq/c,aq/d,aq/e,q/b,q/c,q/d,q/e,qa^2/bcde;q)_{\infty}}.
\end{eqnarray*}

Many proofs of above identity have been found, see, for example,
Slater and Lakin \cite{SlLa}, Andrews \cite{Andrews},  Askey and
Ismail \cite{AsIs}, Askey \cite{Askey}, Chen and Liu
\cite{ChenLiu}, Schlosser \cite{SCHL} and Jouhet and Schlosser
\cite{SCHL04}. Our proof shows that the semi-finite form of the
${}_6\psi_6$ summation is in essence a shifted version of the
$_8\phi_7$ summation. This proof utilizes Ramanujan's ${}_1\psi_1$
summation (\ref{Ran}). It would be interesting to find a proof
using a semi-finite (or even finite) form which yields Bailey's
${}_6\psi_6$ summation in a direct limit, without the need to
invoke another summation formula as above.

\vspace{.33cm}
 \noindent{\bf Acknowledgments.} We thank R. Askey,
Wenchang Chu and M. E. H. Ismail for their valuable comments. In
particular, we are indebted to M. Schlosser for crucial
suggestions leading to considerable improvement of an earlier
version. This work was done under the auspices of the 973 Project
on
 Mathematical Mechanization, the Ministry of Education, the Ministry of Science and Technology, and the
National Science Foundation of China.

\end{document}